\documentclass[11pt]{article}
\usepackage{amsmath, amsthm, amssymb, eucal, setspace, mathrsfs}

\makeatletter
\renewcommand\section{\@startsection{section}{1}{\z@}%
 						{-3.5ex \@plus -1ex \@minus -.2ex}
						{2ex \@plus.2ex}
						{\large\bfseries}}
\renewcommand\subsection{\@ifstar
						{\setcounter{subsection}{\value{equation}}
					\@startsection{subsection}{2}{\z@}
                          {1.75ex \@plus.5ex \@minus.2ex}%
                           {-.4em}		
					\textit*}
					{\setcounter{subsection}{\value{equation}}
						\stepcounter{equation}
					\@startsection{subsection}{2}{\z@}
                          {1.75ex \@plus.5ex \@minus.2ex}%
                           {-.4em}		
					\textit}}
\def\@seccntformat#1{\@ifundefined{#1@cntformat}%
	{\csname the#1\endcsname\quad} 
	{\csname #1@cntformat\endcsname}} 
\def\section@cntformat{\thesection.~} 
\def\subsection@cntformat{(\thesubsection)\ }
\renewcommand*\l@section{\mdseries\small\@dottedtocline{1}{1.5em}{2em}}
\makeatother

\numberwithin{equation}{section}
\theoremstyle{plain}

\swapnumbers
\newtheorem{theorem}[equation]{Theorem}

\theoremstyle{definition}

\theoremstyle{remark}
\newtheorem{remark}[equation]{Remark}

\newcommand{\cA}{\mathscr{A}}

\newcommand{\cP}{\mathscr{P}}

\newcommand{\frg}{\mathfrak{g}}

\newcommand{\frt}{\mathfrak{t}}

\newcommand{\frz}{\mathfrak{z}}

\newcommand{\bC}{\mathbb{C}}

\newcommand{\bR}{\mathbb{R}}
\newcommand{\bT}{\mathbb{T}}
\newcommand{\bZ}{\mathbb{Z}}

\newcommand{\cliff}{\mathrm{Cliff}}
\newcommand{\dirac}{\textup{\mbox{D\hspace{-0.6em}\raisebox{.13ex}{/}\hspace{.07em}}}}
\newcommand{\mfact}{\mathrm{MF}^\tau}
\newcommand{\Pot}{W}
\newcommand{\mi}{\mathrm{i}}

\title{\Large{\textbf{Dirac families for Loop groups as Matrix factorizations}}}

\author{Daniel S.~Freed \and Constantin Teleman}

\begin{document}
\maketitle

\begin{abstract} \noindent {We identify the  category of integrable
lowest-weight representations of the loop group $LG$ of a compact Lie group
$G$ with the linear category of twisted, conjugation-equivariant \emph{curved
Fredholm complexes} on the group $G$: namely, the twisted, equivariant \emph
{matrix factorizations} of a super-potential built from the loop rotation 
action on $LG$. This lifts the isomorphism of $K$-groups of \cite{fht1}--\cite{fht3} 
to an equivalence of categories. The construction uses 
families of Dirac operators.}
\end{abstract}

\section{Introduction and background}
The group $LG$ of smooth loops in a compact Lie group $G$ has a remarkable 
class of linear representations whose structure parallels the theory for 
compact Lie groups \cite{ps}. The defining stipulation is the existence of a circle 
action on the representation, with finite-dimensional eigenspaces and 
spectrum bounded below, intertwining with the loop rotation action 
on $LG$. We denote the rotation circle by $\bT_r$; its infinitesimal 
generator $L_0$ represents the \emph{energy} in a conformal field theory. 

Noteworthy is the \emph{projective nature} of these representations, described 
(when $G$ is semi-simple) by a \emph{level} $h\in H^3_G(G;\bZ)$ in the equivariant 
cohomology for the adjoint action of $G$ on itself. The representation 
category $\mathfrak{Rep}^h(LG)$ at a given level $h$ is semi-simple, with finitely 
many irreducible isomorphism classes. Irreducibles are classified by their 
\emph{lowest weight} (plus some supplementary data, when $G$ is not simply 
connected \cite[Ch.~IV]{fht3}). 

In a series of papers \cite{fht1}--\cite{fht3}\, the authors, jointly with
Michael Hopkins, construct $K^0\mathfrak{Rep}^h(LG)$ in terms of
a twisted, conjugation-equivariant topological $K$-theory group. To wit, when
$G$ is connected, as we shall assume throughout this paper,\footnote {Twisted
loop groups show up when $G$ is disconnected \cite{fht3}.} we have
\begin{equation}\label{oldtheorem} K^0\,\mathfrak{Rep}^h(LG) \cong 
K^{\tau +\dim G}_G(G),
\end{equation}
with a twisting $\tau\in H^3_G(G;\bZ)$ related to $h$, as explained below.
\begin{remark} One loop group novelty is a \emph{braided tensor}
structure\footnote{When $G$ is not simply connected, there is a constraint on
$h$.} on $\mathfrak{Rep}^h(LG)$.  The structure arises from the \emph{fusion
product} of representations, relevant to $2$-dimensional conformal field
theory. The $K$-group in \eqref{oldtheorem} carries a Pontryagin product, and
the multiplications match in \eqref{oldtheorem}.
\end{remark}

The map from representations to topological $K$-classes is implemented by the 
following \emph{Dirac family}. Calling $\cA$ the space of connections on
the trvial $G$-bundle over $S^1$, the quotient stack $[G\!:\!G]$ under 
conjugation is equivalent to $[\cA\!:\!LG]$ under the gauge action, via the 
holonomy map $\cA\to G$. Denote by $\mathbf{S}^\pm$ the (lowest-weight) modules 
of spinors for the loop space $L\frg$ of the Lie algebra and by $\psi(A): 
\mathbf{S}^\pm\to \mathbf{S}^\mp$ the action of a Clifford generator $A$, 
for $d+Adt\in \cA$. A representation $\mathbf{H}$ of $LG$ leads to a family 
of Fredholm operators over $\cA$,
\begin{equation}\label{diracfam}
\dirac_A: \mathbf{H}\otimes\mathbf{S}^+ \to \mathbf{H}\otimes\mathbf{S}^-, \quad
	\dirac_{A} :=\dirac_0 + \mi\psi(A)
\end{equation}
where $\dirac_0$ is built from a certain Dirac operator \cite{land} on the loop 
group.\footnote{It is also the square root $G_0$, in the super-Virasoro algebra, of 
the infinitesimal circle generator $L_0$.} The family is projectively $LG$-equivariant; 
dividing out by the subgroup $\Omega G\subset LG$ of based loops leads 
to a projective, $G$-equivariant Fredholm complex on $G$, whose $K$-theory class 
$\left[(\dirac_\bullet, \mathbf{H}\otimes\mathbf{S}^\pm)\right]\in  K^{\tau +*}_G(G)$ 
is the image of $\mathbf{H}$ in the isomorphism \eqref{oldtheorem}. When $\dim G$ is odd, 
$\mathbf{S}^+=\mathbf{S}^-$ and skew-adjointness of $\dirac_A$ leads to 
a class in $K^1$. The twisting $\tau$ is the level of $\mathbf{H}\otimes\mathbf{S}$ 
as an $LG$-representation, with a ($G$-dependent) shift from the level $h$ 
of $\mathbf{H}$.

The degree-shift is best explained in the world of super-categories, 
with $\bZ/2$ gradings on morphisms and objects; odd simple objects have as
endomorphisms the rank one Clifford algebra $\cliff(1)$, and contribute, in 
the semi-simple case, a free generator to $K^1$ instead of $K^0$.
Consider the $\tau$-projective representations of $LG$ with compatible action 
of $\cliff(L\frg)$, thinking of them as modules for the (not so well-defined) 
crossed product $LG\ltimes\cliff(L\frg)$. They form a semi-simple super-category 
$\mathfrak{SRep}^\tau$, and the isomorphism \eqref
{oldtheorem} becomes 
\begin{equation}\label{oldsupertheorem}
K^*\,\mathfrak{SRep}^\tau\left(LG\ltimes\cliff(L\frg)\right) \cong  K^{\tau +*}_G(G)
\end{equation}
with the advantage of no shift in degree or twisting.\footnote{For simply connected 
$G$, both sides live in degree $\dim\frg$; both parities can be present for general 
$G$.}  This isomorphism is induced 
by the Dirac families of \eqref{diracfam}: a super-representation $\mathbf{SH}^\pm$ 
of $LG\ltimes\cliff(L\frg)$ can be coupled to the Dirac operators $\dirac_\bullet$ 
without a choice of factorization $\mathbf{H} \otimes\mathbf{S}^\pm$.

\subsection*{Acknowledgements.} This research was partially supported by NSF-funded 
Focussed Research Group grants DMS-1160461 and DMS-1160328.

\section{The main result}
There is a curious lack of symmetry in \eqref{oldsupertheorem}: 
the isomorphism is induced by a functor of underlying Abelian categories, 
from $\bZ/2$-graded representations to twisted Fredholm bundles over $G$, 
but this functor is far from an equivalence. The category $\mathfrak{SRep}^\tau$ 
is semi-simple (in the graded sense discussed), but that of twisted Fredholm bundles 
is not so. We can even produce continua of non-isomorphic objects 
in any given $K$-class by perturbing a Fredholm family 
by compact operators. 

Here, we redress this problem via the inclusion of a \emph{super-potential 
$\Pot$}, of algebraic geometry and $2$-dimensional physics $B$-model fame. As 
explained by Orlov\footnote{Orlov discusses complex algebraic vector bundles; we 
found no suitable exposition covering equivariant Fredholm complexes in topology, 
and a discussion is planned for our follow-up paper.} \cite{orlov}, this 
deforms the category of complexes of vector bundles into that of \emph{matrix 
factorizations}: the \emph{$2$-periodic, curved complexes} with curvature $\Pot$. 
Our super-potential will have Morse critical points, leading to a semi-simple 
super-category with one generator for each critical point. The generators 
are precisely the Dirac families of \eqref{diracfam} on  irreducible 
$LG$-representations. This artifice of introducing a super-potential is redeemed 
by its natural topological origin, in the \emph{loop rotation} action on 
the stack $[G\!:\!G]$. (Rotation is more evident in the presentation by connections, 
$[\cA\!:\! LG]$.) Namely, for twistings $\tau$ transgressed from $BG$, the action of 
$\bT_r$  refines to a $B\bZ$-action on the $G$-equivariant \emph{gerbe} $G^\tau$ 
over $G$ which defines the $K$-theory twisting. The logarithm of this lift is 
$2\pi\mi W$. 
\begin{remark}
\begin{enumerate}
\item The conceptual description of a super-potential as logarithm of a $B\bZ$-action 
on a category of sheaves is worked out in \cite{p}; the matrix factorization 
category is the \emph{Tate fixed-point category} for the $B\bZ$-action. 
On varieties, $W$ is a function and $\exp(2\pi\mi W)$ defines the 
$B\bZ$-action; on a stack, there can be (as here) a geometric underlying 
action as well.
\item To reconcile our story with \cite{p}, we must  rescale $W_\tau$ 
so that it takes integer values at all critical points; we will ignore this 
detail to better connect with the formulas in \cite{fht2, fht3}. 
\end{enumerate}
\end{remark}
To spell out our construction, recall that a stack is an instance of a
category. A $B\bZ$-action thereon is described by its generator, an
automorphism of the identity functor. This is a section over the space of
objects, valued in automorphisms, which is central for the groupoid
multiplication. For $[G\!:\!G]$, the relevant section is the identity
map $G\to G$ from objects to morphisms. Intrinsically, $[G\!:\! G]$ is the
mapping stack from $B\bZ$ to $BG$, and the $B\bZ$-action in question is the
self-translation action of $B\bZ$.  This rigidifies the $\bT_r$-action on the
homotopy equivalent spaces $LBG\sim BLG\sim \cA/LG$.

A class $\hat{\tau}\in H^4(BG;\bZ)$ transgresses to a $\tau\in H^3_G(G;\bZ)$, 
with a natural $\bT_r$-equivariant refinement. This can also be rigidified, 
as follows. The exponential sequence lifts $\hat{\tau}$ uniquely to 
$H^3(BG; \bT)$, the group cohomology with smooth circle coefficients. That  
defines a Lie $2$-group $G^{\hat{\tau}}$, a multiplicative $\bT$-gerbe 
over $G$. (Multiplicativity encodes the original $\hat{\tau}$.) 
The mapping stack from $B\bZ$ to $BG^{\hat{\tau}}$ is the quotient  
$[G^{\hat{\tau}}\!:\!G^{\hat{\tau}}]$ under conjugation, and carries a natural 
$B\bZ$-action from the self-translations of the latter. Because 
$B\bT\hookrightarrow G^{\hat{\tau}}$ is strictly central, the self-conjugation 
action of $G^{\hat{\tau}}$ factors through $G$, and the quotient stack 
$[G^{\hat\tau}\!:\!G]$ is our $B\bZ$-equivariant gerbe over $[G\!:\! G]$ with band $\bT$. 
We denote this central circle by $\bT_c$, to distinguish it from $\bT_r$. 

The $B\bZ$-action gives an automorphism $\exp(2\pi\mi\Pot_\tau)$ of the 
identity of $[G^{\hat\tau} \!:\!G]$, lifting the one on $[G\!:\!G ]$. Concretely, 
$[G^{\hat\tau} \!:\!G]$ defines a $\bT_c$-central extension of the stabilizer 
of $[G\!:\! G]$, and $\exp(2\pi\mi\Pot_\tau)$ is a trivialization of its 
fiber over the automorphism $g$ at the point $g\in G$; see \S3 below. 
The logarithm $\Pot_\tau$ is multi-valued and only locally well-defined; nevertheless, 
the category of twisted matrix factorizations, $\mfact_G(G; W_\tau)$ is well-defined, 
and its objects are represented by $\tau$-twisted $G$-equivariant Fredholm complexes 
over $G$ curved by $\Pot_\tau+\bZ\cdot\mathrm{Id}$. 

\begin{theorem}\label{main}
The following defines an equivalence of categories from $\mathfrak{SRep}^\tau$ 
to $\mfact_G(G; W_\tau)$: 
a graded representation $\mathbf{SH}^\pm$ goes to the twisted and curved 
Fredholm family $\left(\dirac_\bullet,\mathbf{SH}^\pm\right)$ whose value at 
the connection $d+A\,dt\in \cA$ is the $\tau $-projective $LG$-equivariant curved Fredholm complex
\[
\dirac_{A} =\dirac_0 + \mi\psi(A): \mathbf{SH}^+\rightleftarrows \mathbf{SH}^- 
\]
\end{theorem}

\begin{remark}
\begin{enumerate}
\item Matrix factorizations obtained from irreducible representations are 
supported on single conjugacy classes, the so-called \emph{Verlinde conjugacy 
classes} in $G$, for the twisting $\tau$. These are the supports of the co-kernels 
of the Dirac families \eqref{diracfam}, \cite[\S12]{fht3}.
\item There is a braided tensor structure on $\mathfrak{SRep}^\tau \left(LG\ltimes
\cliff(L\frg)\right)$ (without $\bT_r$-action). A corresponding structure 
on $\mfact_G(G,\Pot_\tau)$ should come from the Pontryagin product. We do 
not know how to spell out this structure, partly because the $\bT_r$-action is 
already built into the construction of $\mfact$, and the Pontryagin 
product is \emph{not} equivariant thereunder.  
\item The values of the automorphism $\exp(2\pi\mi W_\tau)$ at the Verlinde conjugacy 
classes determine the \emph{ribbon element} in $\mathfrak{Rep}^h(LG)$; see \cite{fhlt} for 
the discussion when $G$ is a torus.
\end{enumerate}
\end{remark}

Theorem~\ref{main} has a $\hat{\tau}\to\infty$ scaling limit, which we will use 
in the proof. In this limit, the representation category of $LG$ becomes that 
of $G$. On the topological side, noting that each $\hat{\tau}$ defines an 
inner product  on $\frg$, we zoom into a neighborhood of $1\in G$ so that the 
inner product stays fixed. This leads to a $G$-equivariant matrix 
factorization category $\mathrm{MF}_G(\frg,\Pot)$ on the Lie algebra. The 
$\tau$-central extensions of stabilizers 
near $1$ have natural splittings, and the $W_\tau$ converge to a super-potential 
$\Pot\in G\ltimes\mathrm{Sym}(\frg^*)$, which, in a basis $\xi_a$ of $\frg$ with 
dual basis $\xi^a$ of $\frg^*$, we will calculate in~\S3 to be
\begin{equation}\label{liepotential}
\Pot = -\mi\cdot\xi_a(\delta_1)\otimes \xi^a + \frac{1}{2} \sum\nolimits_a\|\xi^a\|^2
\end{equation} 
with $\xi_a(\delta_1)$ denoting the respective derivative of the delta-function 
at $1\in G$. It is important that \eqref{liepotential} is central 
in the crossed product algebra $G\ltimes\mathrm{Sym}(\frg^*)$.

To describe the limiting case, recall from \cite[\S4]{fht3} the $G$-analogue of the Dirac 
family \eqref{diracfam}. Kostant's \emph{cubic Dirac operator} \cite{kos} on $G$ 
is left-invariant, and the Peter-Weyl decomposition gives 
an operator $\dirac_0: \mathbf{V}\otimes \mathbf{S}^\pm\to \mathbf{V}\otimes \mathbf{S}^\mp$ 
for any irreducible representation $\mathbf{V}$ of $G$, coupled to the spinors 
$\mathbf{S}^\pm$ on $\frg$. As before, it is better to work with graded 
modules for the super-algebra $G\ltimes\cliff(\frg)$. 

\begin{theorem}\label{limit}
Sending $\mathbf{SV}^\pm$ to $\left(\dirac_\bullet,\mathbf{SV}^\pm\right)$, 
the curved complex over $\frg$ given by
\[
\frg\ni\mu\mapsto \dirac_\mu = \dirac_0 + \mi\psi(\mu):\mathbf{SV}^+\leftrightarrows \mathbf{SV}^-
\] 
provides an equivalence of super-categories from graded $G\ltimes\cliff(\frg)$-modules 
$\mathbf{SV}^\pm$ to $G$-equivariant, $\Pot$-matrix factorizations over $\frg$.
\end{theorem}
\noindent
With  $\lambda$ denoting the lowest weight of $V$ and $T(\mu)$ the $\mu$-action on 
$\mathbf{SV}$, we have  \cite[Cor.~4.8]{fht3}
\[
\dirac_\mu^2 = -\|\lambda_V+\rho\|^2 +2\mi\cdot T(\mu) -\|\mu\|^2 \in (-2W) +\bZ.
\]

\section{Outline of the proof}

\subsection{Executive summary.} 
The category $\mfact_G(G;W_\tau)$ sheafifies over the conjugacy classes of $G$. 
Near any $g\in G$ with centralizer $Z$, the stack $[G\!:\!G]$ is modeled on a 
neighborhood of $0$ in the adjoint quotient $[\frz\!:\!Z]$ of the Lie algebra $\frz$, 
via $\zeta\in \frz\mapsto g\cdot\exp(2\pi\zeta)$. We will compute the local 
$W_\tau$ in the crossed product $Z\ltimes C^\infty\left({\frz}\right)$, 
recovering \eqref{liepotential}, up to a $g$-dependent central translation 
in $\frz$. We then show that $\mfact$ lives only on \emph{regular} elements 
$g$. Assuming for brevity that $\pi_1(G)$ is torsion-free: $Z$ is then the 
maximal torus $T\subset G$, where the super-potential $\Pot_\tau$ turns 
out to have Morse critical points, located precisely at the Verlinde conjugacy 
classes. The local category is freely 
generated by the respective Atiyah-Bott-Schapiro Thom complex.\footnote
{The Clifford multiplication  acts in both directions, giving a curved complex.} 
The latter is quasi-isomorphic to our Dirac family for a specific 
irreducible representation, associated with the Verlinde class  \cite[\S12]{fht3}.
   
\subsection{The $2$-group.} 
We  use a  \emph{Whitehead crossed module} \cite{white} description for 
$G^{\hat{\tau}}$. This is an exact sequence of groups
\[
\bT_c\rightarrowtail K\xrightarrow{\:\varphi\:} H \twoheadrightarrow G,  
\]
equipped with an action of $H$ on $K$ which lifts the self-conjugation of $H$ and 
factors the self-conjugation of $K$ via $\varphi$. Call $h$ an $H$-lift of $g$, 
and $C$ the pre-image of $Z$ in $H$. Define the central extension $\widetilde
Z$ 
by means of a $\bT_c$-central extension of $C$, trivialized over $\varphi(K)\cap C$, 
as follows. The commutator $c\mapsto hc h^{-1}c^{-1}$ gives crossed homomorphism 
$\chi: C \to \varphi(K)$, with respect to the conjugation action of $C$ on 
$\varphi(K)$. The action having been lifted to $K$, $\chi$ pulls back the 
central extension $\bT_c\rightarrowtail K\to \varphi(K)$ to $C$. The $h$-action 
on $K$ identifies the fiber of $K$ over any $c\in \varphi(K)$ with that over 
$hc h^{-1}$, trivializing the pull-back extension over $\varphi(K)$. Finally, 
$hhh^{-1}h^{-1}=1$, so the extension is also trivialized over $c=h$, 
defining our $\exp(2\pi\mi W_\tau)$ at $g\in Z$.

\subsection{Computing the local super-potential.} 
Following~\cite{bscs},
take  $K=\Omega^\tau G$, the $\tau$-central extension of the group of based 
smooth maps $[0,2\pi]\to G$ sending $\{0,2\pi\}$ to $1$, and $H=\cP_1G$, the group 
of smooth paths $[0,2\pi]\to G$ starting at $1\in G$. The requisite 
$H$-action on the Lie algebra $\mi\bR\oplus\Omega \frg$ of $K$ is 
\begin{equation}\label{cocycle}
\gamma.(x\oplus \alpha) = \left(x-\frac{\mi}{2\pi}\int_0^{2\pi}\langle \gamma^{-1}d\gamma |
\alpha\rangle \,\oplus \, \mathrm{Ad}_\gamma (\alpha)\right)
\end{equation} 
 extending the Ad-action of $\Omega^\tau G$ \cite[Prop.~4.3.2]{ps}, and 
exponentiating to an action on $\Omega^\tau G$. (Acting on other components 
of $\Omega G$ requires topological information from $\hat{\tau}$.)

The equivariant gerbe $[G^{\hat\tau}\!:\!G]$ is locally trivialized (possibly on a 
finite cover of $Z$) uniquely up to discrete choices: 
the automorphisms of the central extension $\widetilde{Z}$. 
We spell out $W_\tau$ in these terms. Lift $g$ to $\cP_1G$ as $h=\exp(t\mu)$, 
for a shortest logarithm $2\pi\mu$ of $g$, and assume for now that $Z$ 
centralizes $\mu$. Instead of the entire group $C$, we use in the construction 
the subgroup $\cP_1Z$ of paths in $Z$. It centralizes $h$, and this trivializes
our $\bT_c$-extension over $\cP_1Z$, with $W_\tau=0$. However, by \eqref{cocycle}, 
the extension over $\Omega Z = \varphi(K)\cap\cP_1Z$ is trivialized by 
the Lie algebra character $\alpha\mapsto -\frac{\mi}{2\pi}\int_0^{2\pi} \langle 
\mu|\alpha\rangle dt$. To trivialize $\widetilde{Z}$, we must therefore extend this 
to a linear character of $\cP^{\phantom M}_0\!\!\frz$. The same formula \eqref{cocycle} does this, 
supplying the locally constant trivialization of $\widetilde{Z}$. We now get the value 
$2\pi\mi W_\tau(g) = \pi\mi\|\mu\|^2\oplus2\pi\mu \in \mi\bR\oplus\frg$. 

At the remaining points, $W_\tau$ is determined by continuity, but 
can also be pinned down by the restriction to a maximal torus containing $g$.

\subsection{Vanishing of singular contributions.}\label{vanish}
When $\frz$ is non-abelian, we show the vanishing of the matrix factorization 
category localized at $g$. Take $g=1, Z=G$, $W$ on $\frg$ as in \eqref{liepotential}, 
plus possibly a central linear term $\mu$. Koszul duality equates the localized 
category $\mfact_G(\frg;W)$ with the super-category of $\bZ/2$-graded modules over 
the differential super-algebra $(G\ltimes\cliff(\frg), [\dirac_\mu,\,\underline{\ }\,])$; 
$\dirac_\mu = \dirac_0+\mi\psi(\mu)$, with Kostant's cubic Dirac operator of $\S2$. 
Ignoring the differential, the algebra is semi-simple, with simple 
modules the $\mathbf{V}\otimes\mathbf{S}^\pm$ of Theorem~\ref{limit}, 
for the irreducible $G$-representations $\mathbf{V}$. Since 
$\dirac_\mu^2 = -\|\lambda_V+\mu+\rho\|^2<0$, $[\dirac_\mu,\dirac_\mu]$ 
provides a homotopy between $0$ and a central unit in the algebra. 
This makes the super-category of graded modules quasi-equivalent 
to $0$.

\subsection{Globalization for the torus.} We  describe the stack
$[T^{\hat\tau}\!:\!T]$ and potential $W_\tau$ in the presentation $T= [\frt\!:\!\Pi]$ of
the torus as a quotient of its Lie algebra by $\Pi\cong\pi_1(T)$. Lifted to~
$\frt$, the gerbe of stabilizers $\tilde{T}$ is trivial, with band $T\times\bT_c$.  The
descent datum under translation by $p\in\Pi$ is the shearing automorphism of
$T\times \bT_c$ given by the character $t\mapsto\exp\langle p| \log
t\rangle,\;t\in T$. In the same trivialization over $\frt$, the
super-potential is 
\[ 2\pi\mi W_\tau(\mu)= \pi\mi\|\mu\|^2\oplus 2\pi\mu\in
\mi\bR\oplus\frt, 
\] 
the first factor being the Lie algebra of $\bT_c$. That is the function 
``$\frac{1}{2}\|\log\|^2$" on $T$, invariant under $p$-translation, save 
for an additive shift by the integer $\|p\|^2/2$.

With $\Lambda$ denoting the character lattice of $T$, the crossed product
algebra of the stack $[T^\tau:T]$ can be identified with the functions on
$\left(\coprod_{\lambda\in \Lambda} \frt_\lambda\right)/\Pi$, with the action
of $\Pi$ by simultaneous translation on $\Lambda$ and $\frt$. On the sheet
$\lambda\in \Lambda$, $W_\tau = -\langle\lambda | \mu\rangle + \|\mu\|^2/2$
has a single Morse critical point at $\mu = \lambda$.

It follows that the super-category $\mfact_T(T;\Pot_\tau)$ is semi-simple, 
with one generator of parity $\dim\frt$ at each point in the kernel of the 
isogeny $T\to T^*$ derived from the quadratic form $\hat{\tau}\in H^4(BT;\bZ)$. 
The kernel comprises precisely the Verlinde points in $T$ \cite{fhlt},  and
this concludes the proof.

\noindent
\small{\textsc
{D.S.~Freed:} UT Austin, Mathematics Department, RLM 8.100, 2515 Speedway
C1200, Austin, TX 78712 \\
\textsc{C.~Teleman:} UC Berkeley, Mathematics Department, 970 Evans Hall \#3840, 
Berkeley, CA 94720
}

\end{document}